\documentclass[a4paper,12pt]{article}
\usepackage[T2A]{fontenc}
\usepackage[english]{babel}
\usepackage{graphicx}
\usepackage{amsmath}
\usepackage{mathrsfs}
\usepackage{amssymb}
\usepackage{caption}
\usepackage[usenames]{color}
\usepackage{colortbl}

\begin{document}
\begin{center}

\vskip 5mm {\Large \bfseries About one method of parallelization of calculations during the reconstruction of a tomographic image}

\vskip 5mm {\bf A.A. Alikhanov$^1$, A.M. Apekov$^1$, Z.A. Kokov$^{2,3}$, A.O. Belyaev$^{2}$, L.A. Khamukova$^{2,3}$}
\vskip 2mm

$^1${Institute of Applied Mathematics and Automation KBSC RAS, Nalchik, 360000, Russia}
\vskip 2mm
$^{2}${Southern Federal University, Rostov-on-Don, 344006, Russia }
\vskip 2mm
$^{3}${Kabardino-Balkarian State University named after H.M. Berbekov, Nalchik, 360004, Russia}

{\it e-mail: aaalikhanov@gmail.com}

\abstract{A new method for solving systems of linear algebraic equations of a special type arising in solving problems of image reconstruction has been proposed. This method, due to a certain symmetry of the matrix and the choice of the voxel numeration method for two-dimensional problems, allows us to divide the initial system of algebraic equations into two independent systems, which enables us to carry out the calculation in parallel. The size of the matrices of the resulting systems is 4 times less than the size of the original matrix, and these matrices are less dispersed.}


\end{center}
\section{Introduction}

An important method of non-invasive imaging widely used in medical diagnostics is x-ray computed tomography. Image reconstruction is performed both by direct and iterative methods. 
Iterative methods are preferred when using a limited set of projections, for example, in digital tomosynthesis.
Digital tomosynthesis is a new method of three-dimensional reconstruction of an object using projections at limited viewing angles. 
Currently, the method of simultaneous algebraic reconstruction  (SART)  is a generally accepted iterative reconstruction method \cite{001}. 
The disadvantage of SART is its high computational complexity and cost. 
In the works \cite{002,003,004}, the possibility of improving the SART image reconstruction speed was investigated. 
The use of parallel methods to speed up the SART reconstruction was proposed in \cite{005,006,007,008,009,010, ZhangSL}.

In \cite{alikh} the method of partitioning the reconstruction region into an uneven voxel grid was considered. This allows one to split the original three-dimensional problem into a chain of two-dimensional independent tasks.

{In this paper, we propose a method that, by virtue of a certain symmetry of the matrix and the choice of the numbering method for voxels for two-dimensional problems, allows us to split the initial system of algebraic equations into  two  independent systems, which makes it possible to perform a parallel calculation. The size of the matrices of the resulting systems is 4 times less than the size of the original matrix, and these matrices are less dispersed.}

\section{Methods}

Let $f = (f_1, \ldots, f_N) \in R^{N}$ be a discrete representation of an unknown image ($N$ is the total number of image voxels) to be reconstructed, $p = (p_1, \ldots, p_M) \in R^{M}$ be the measured projection data ($M$ is the total number of rays), and $A = (w_{i, j})$ be a known system matrix, whose element $w_{i,j}$ is the weighting factor that represents the contribution of the $j$th voxel of the $i$th ray integral. Hereinafter we shall assume that the numbers $N$ and $M$ are even (as a rule, it takes place in practice).
The image reconstruction problem can be formulated as a system of linear equations:
\begin{equation}
\label{ur1}
A f = p
\end{equation}

Suppose that during the motion of the radiator, the irradiation is performed in positions that are pairwise symmetrical about the axis OZ. let the number of radiation is  $K$.
The reconstruction region is divided into layer-by-layer voxels, and each of the $j$ voxels is characterized by the parameter $f_j$.
However, voxels are numbered in a special order (according to the vertical "snake" rule) as shown in Fig. 1 for the case $N=32$.

\begin{figure}[h]
\center{\includegraphics[scale=0.24]{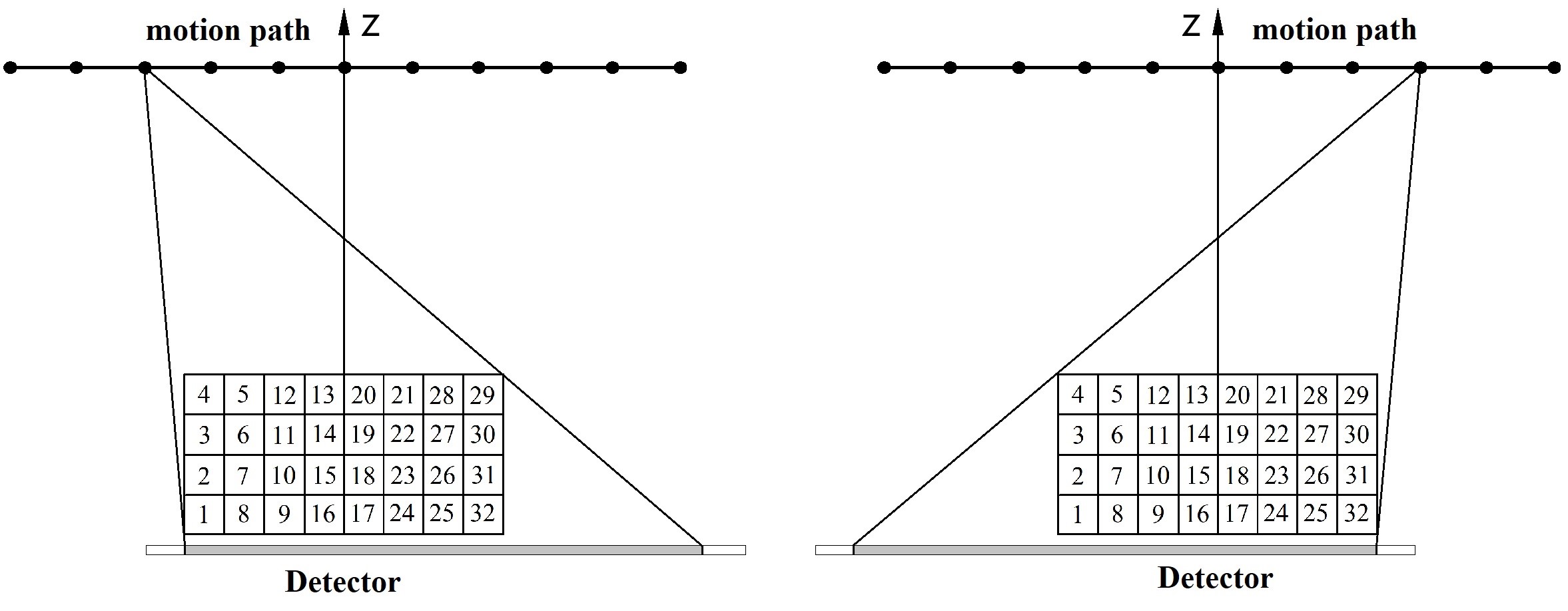}}
\caption{Rule of voxel numbering in the reconstruction plane. }
\label{fig1}
\end{figure}

In this case, the matrix $А$ is symmetric in the sense
\begin{equation}
\label{ur2}
w_{i, j} = w_{M-i+1, N-j+1}.
\end{equation}

Consider two systems of linear algebraic equations that are independent of each other:

\begin{equation}
\label{ur3}
A^{(1)} f^{(1)} = p^{(1)},
\end{equation}

\begin{equation}
\label{ur4}
A^{(2)} f^{(2)} = p^{(2)},
\end{equation}
where the elements of matrices $A^{(1)}$ and $A^{(2)}$ depend on the elements of the matrix A, and are found by the formulas:
\begin{equation}
\nonumber
w_{i,j}^{(1)} = w_{i,j} - w_{M-i+1,j}, \quad w_{i,j}^{(2)} = w_{i,j} + w_{M-i+1,j},\quad i = 1,\ldots, \frac{M}{2},\quad j = 1,\ldots, \frac{N}{2},
\end{equation}
and the elements of vectors $p^{(1)}$ and $p^{(2)}$ depend on the elements of vector $p$:
\begin{equation}
\nonumber
p_{i}^{(1)} = {p_{i} - p_{M-i+1}}, \quad p_{i}^{(2)} = {p_{i} + p_{M-i+1}},\quad i = 1,\ldots, \frac{M}{2}. 
\end{equation}

\textbf{Theorem 1.} \textit{Let the matrix of system (\ref{ur1}) satisfy the symmetry conditions (\ref{ur2}), whereas $f^{(1)}$ and $f^{(2)}$ are the solutions of systems (\ref{ur3}) and (\ref{ur4}) respectively, then $f = (f_1, \ldots, f_N)$,  where
 \begin{equation}
 \label{ur5}
 f_j = \frac{f_j^{(2)} + f_j^{(1)}}{2},\quad
 f_{N-j+1} = \frac{f_j^{(2)} - f_j^{(1)}}{2}, \quad j = 1,\ldots, \frac{N}{2},
 \end{equation}
is a solution to equation (\ref{ur1}).}

\textbf{Proof.} We introduce the notation
$$
f_j^{(1)} = {f_j - f_{N-j+1}}, \quad f_j^{(2)} = {f_j + f_{N-j+1}}, \quad j = 1,\ldots, {N}.
$$
It's obvious that 
\begin{equation}
\label{ur6}
f_j^{(1)} = - f_{N-j+1}^{(1)}, \quad f_j^{(2)} = f_{N-j+1}^{(2)}, \quad j = 1,\ldots, {N},
\end{equation}
From equation (\ref{ur1}), by virtue of the equalities $f = \frac{f^{(1)} + f^{(2)}}{2}$ and $p = \frac{p^{(1)} + p^{(2)}}{2}$, we get
\begin{equation}
\label{ur7}
A \left(\frac{f^{(1)} + f^{(2)}}{2}\right) = \frac{1}{2}A f^{(1)} + \frac{1}{2}A f^{(2)} = \frac{1}{2}p^{(1)} + \frac{1}{2}p^{(2)},
\end{equation}
where $p^{(1)}$ and $p^{(2)}$ are extended for all $i = 1,\ldots, {N}$ so that the symmetry conditions are satisfied (\ref{ur6}).

By virtue of the symmetry condition (\ref{ur2}) for matrix $A$, vector $A f^{(1)}$ satisfies the first symmetry condition from (\ref{ur6}), whereas vector $A f^{(2)}$  satisfies the second symmetry condition from (\ref{ur6}). Since the vector representation as the sum of two vectors satisfying the symmetry conditions (\ref{ur6}) is unique, then (\ref{ur7}) splits into two systems:
\begin{equation}
\label{ur8}
A f^{(1)} = p^{(1)}, 
\end{equation}
\begin{equation}
\label{ur9}
A f^{(2)} = p^{(2)}. 
\end{equation}

System (\ref{ur8}) can be represented as 
\begin{equation}
\label{ur10}
\begin{cases}
A^{(1)} f^{(1)} = p^{(1)},  \\ 
-A^{(1)} f^{(1)} = -p^{(1)}, 
\end{cases} 
\Leftrightarrow\quad
A^{(1)} f^{(1)} = p^{(1)}.
\end{equation}
And system (\ref{ur9}) can be represented as 
\begin{equation}
\label{ur11}
\begin{cases}
A^{(2)} f^{(2)} = p^{(2)},  \\ 
A^{(2)} f^{(2)} = p^{(2)}, 
\end{cases} 
\Leftrightarrow\quad
A^{(2)} f^{(2)} = p^{(2)}.
\end{equation}
Theorem 1 is proved.

The size of the matrices $A^{(1)}$ and $A^{(2)}$ is 4 times less than the size of the original matrix $A$. In addition, the matrices $A^{(1)}$ and $A^{(2)}$ are less dispersed. The independence of these matrices makes it possible to find solutions of systems (3) and (4) by the method of parallel calculations. Obviously, if the number of arithmetic operations required for solving system (\ref{ur1}) is proportional to $N^{\alpha}\cdot M^{\beta}$ ($\alpha > 0,\, \beta > 0$), then the time spent on solving system (\ref{ur1}) using the same method with the help of (\ref{ur5}) will be $2^{\alpha + \beta}$ times less.

\textbf{Theorem 2.} \textit{Let $\bar f^{(1)} = (\bar f^{(1)}_{1}, \ldots, \bar f^{(1)}_{\frac{N}{2}})$ - be the normal pseudo-solution of system (\ref{ur3}), and $\bar f^{(2)} = (\bar f^{(2)}_{1}, \ldots, \bar f^{(2)}_{\frac{N}{2}})$ be the normal pseudo-solution of system (\ref{ur4}), then 
\begin{equation}
\label{ur11_1}
\bar f_j = \frac{\bar f_j^{(2)} + \bar f_j^{(1)}}{2},\quad
\bar f_{N-j+1} = \frac{\bar f_j^{(2)} - \bar f_j^{(1)}}{2}, \quad j = 1,\ldots, \frac{N}{2},
\end{equation}
is a normal pseudo-solution of equation (\ref{ur1}).}

The proof of Theorem 2 follows from the following obvious equalities:
$$
\|f\|_2^2 = \sum\limits_{j=1}^{N}f_j^2 = \sum\limits_{j=1}^{N/2}\left(\left(\frac{f_j^{(2)} + f_j^{(1)}}{2}\right)^2 + \left(\frac{f_j^{(2)} - f_j^{(1)}}{2}\right)^2\right) =
$$
$$
\frac{1}{2}\sum\limits_{j=1}^{N/2}(f_j^{(1)})^2 + \frac{1}{2}\sum\limits_{j=1}^{N/2}(f_j^{(2)})^2 =
\frac{1}{2} \|f^{(1)}\|_2^2 + \frac{1}{2} \|f^{(2)}\|_2^2.
$$

\textbf{Example 1.} Consider the case $N = 6$, $M = 4$
\begin{equation}
\label{ur12}
\left(\begin{array}{cccccc}
	1 & 3 & 5 & 7 & 9 & 1\\
	2 & 4 & 6 & 8 & 3 & 7\\
	7 & 3 & 8 & 6 & 4 & 2\\
	1 & 9 & 7 & 5 & 3 & 1
\end{array}\right)
\left(\begin{array}{c}
f_{1} \\
f_{2} \\
f_{3} \\
f_{4} \\
f_{5} \\
f_{6} 
\end{array}\right)
=
\left(\begin{array}{c}
5 \\
6 \\
8 \\
7 
\end{array}\right).
\end{equation}

In order to solve system (\ref{ur12}) consider the following two systems:
\begin{equation}
\label{ur12_1}
\left(\begin{array}{ccc}
0 & -6 & -2 \\
-5 & 1 & -2
\end{array}\right)
\left(\begin{array}{c}
f_{1}^{(1)} \\
f_{2}^{(1)} \\
f_{3}^{(1)}  
\end{array}\right)
=
\left(\begin{array}{c}
-2 \\
-2 
\end{array}\right),
\end{equation}

\begin{equation}
\label{ur12_2}
\left(\begin{array}{ccc}
2 & 12 & 12 \\
9 & 7 & 14
\end{array}\right)
\left(\begin{array}{c}
f_{1}^{(2)} \\
f_{2}^{(2)} \\
f_{3}^{(2)}  
\end{array}\right)
=
\left(\begin{array}{c}
12 \\
14 
\end{array}\right).
\end{equation}

The general solution of systems (\ref{ur12_1}) and (\ref{ur12_2}) respectively has the form:
\begin{equation}
\label{ur12_1_1}
f_1^{(1)} = 14x, \quad f_2^{(1)} = 10x, \quad f_3^{(1)} = -30x+1,
\quad \text{where} \quad x\in R,
\end{equation}
\begin{equation}
\label{ur12_2_2}
f_1^{(2)} = 84y, \quad f_2^{(2)} = 80y, \quad f_3^{(2)} = -94y+1,
\quad \text{where} \quad y\in R.
\end{equation}
 From (\ref{ur12_1_1}) and (\ref{ur12_2_2}) based on Theorem 1, we obtain the general solution of system (\ref{ur12})
 $$
 f=
 \left(\begin{array}{c}
 7x + 42y             \\
 5x + 40y             \\
 -15x -47y +1          \\
 15x -47y         \\
 -5x + 40y            \\
 -7x+42y
 \end{array}\right).
 $$
Normal pseudo-solutions, up to four decimal places, of systems (\ref{ur12_1}) and (\ref{ur12_2}) respectively, have the form:
$$
f^{(1)} = \left(0.3512, 0.2508, 0.2475\right)^T, \quad \|f^{(1)}\|_2 = 0.4975,
$$
$$
f^{(2)} = \left(0.3542, 0.3373, 0.6036\right)^T, \quad \|f^{(2)}\|_2 = 0.7769,
$$
Normal pseudo solutions, up to four decimal places, for the initial system (\ref{ur12}) have the form:
$$
f = \left(0.3527, 0.2941, 0.4256, 0.1781, 0.0433, 0.0015\right)^T, \quad \|f\|_2 = 0.6523.
$$
For this example, the validity of Theorem 2 is confirmed.

Numerous calculations for square matrices show the validity of the following equality 
\begin{equation}
\label{ur13}
\det(A) = \det(A^{(1)})\cdot\det(A^{(2)}).
\end{equation}

If the matrix $A$ is not square and satisfies symmetry conditions in the sense of (\ref{ur2}), then the matrix
$B = A^{T}\cdot A$ will also be symmetric in the sense of (\ref{ur2}). In addition, the following equalities hold true:
$$
B^{(1)} = (A^{(1)})^T \cdot A^{(1)}, \quad B^{(2)} = (A^{(2)})^T \cdot A^{(2)}.
$$

\textbf{Example 2.} The emitter rotates about a fixed axis parallel to the detector matrix. An object is fixed between the emitter and detector matrix. The lower boundary of the object is at the distance of $h_e=1$ m from the emitter and at a distance of $h_m=0.25$ m from the detector matrix. The size of the matrix is $L_m\times L_m$
($L_m=0.43$ m) and the resolution is $N_p\times N_p$ ($N_p=1024$ pixels). The spot size of an x-ray tube is $a_e=7\cdot 10^{-4}$ m. When rotating the emitter the angle of inclination changes from $-\gamma$ to $+\gamma$, where $\gamma=30^\circ$ (see Fig. 2).

\begin{figure}[h]
	\center{\includegraphics[scale=0.35]{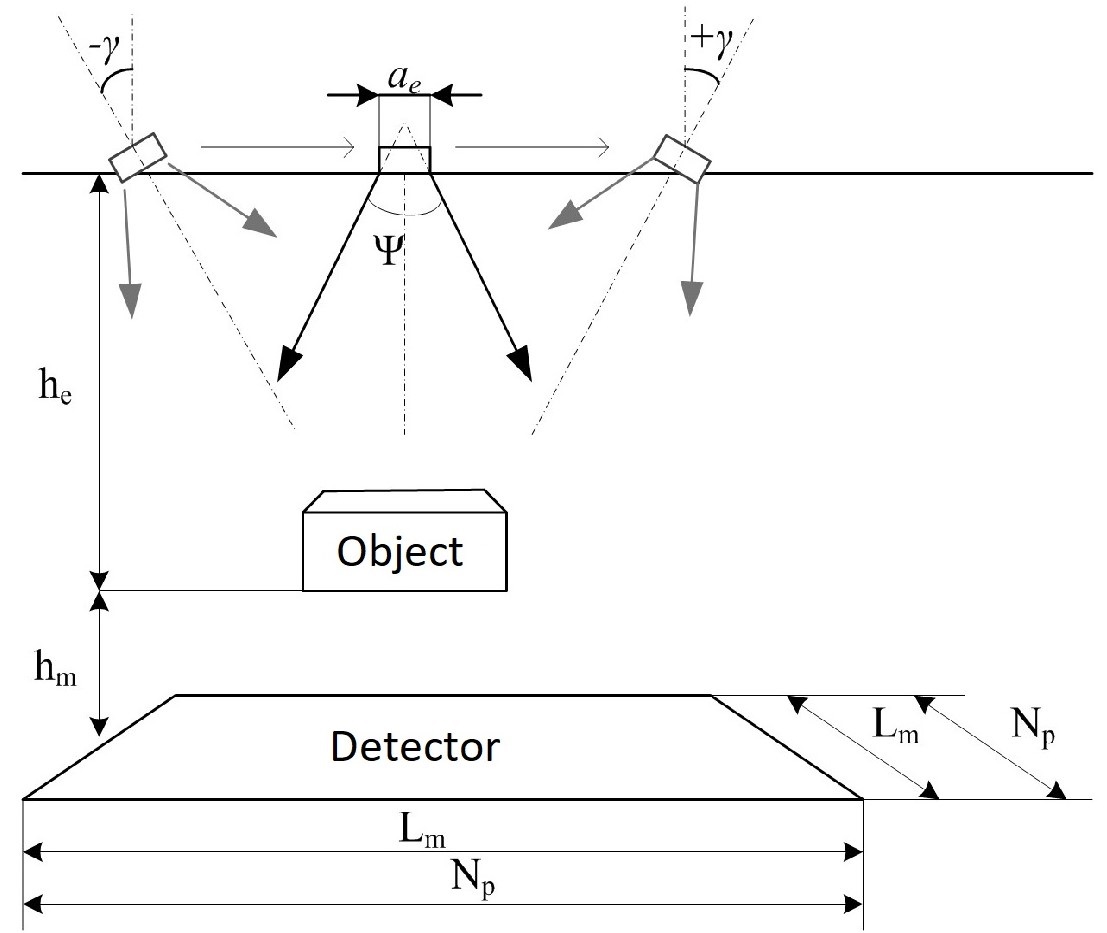}}
	\caption{The result of the reconstruction of the Shepp-Logan phantom. }
	\label{fig2}
\end{figure}

The passage of the beam through voxels of the object is computationally modeled, thus building the system of equations using system of matrix $A$ and the right-hand side $p$ of the system (1). In this example, the reconstruction of the Shepp-Logan phantom is modeled once the object is non-uniformly divided into a voxel grid, as shown in [12]. This allows to obtain a series of two-dimensional independent problems. In example 2 reconstruction of one of the $2D$ - problem using the method proposed is considered. 
In the case when $N = 32\times32$ for the $2D$ - object restoration by the direct method 12 exposures from points located evenly between the extreme left position and center are sufficient as well as at the points symmetric to them relative to the central axis (24 positions in total). Thus we get matrix A with the size $1024\times2280$.
With each doubling of the number of partitions of the object in each direction, the number of positions doubles, and the size of the matrix $A$ increases 16 times. For example, if $N = 64\times64$ then the exposure is produced from 48 positions and the size of matrix $A$ is $4096\times9120$. While in the case when $N = 128\times128$, the exposure is produced from 96 positions with matrix $A$ size $16384\times36480$.
The calculations obtained for these systems of equations are given in Table 1. 
There are three dimensions for matrix $A$ system, the calculation time of the direct method, the errors of the direct method, the calculation time of the method proposed in this paper, and the errors of the proposed method.
Table 1 shows that the advantages of the proposed method  are in the reduced calculation time by an order of magnitude and also in the calculation accuracy. Fig. 3 shows the reconstruction results of the Schepp-Logan phantom for dividing an object into a) $32\times32$, b) $64\times64$ and c) $128\times128$ voxels. 

\begin{table}
	\caption{CPU time and errors of the calculations}
\begin{tabular}{|c|c|c|c|c|c|}
	\hline
	\textnumero &size of $A$ & CPU time & error & CPU time for         & error for\\
	            &                &          &        & our method   & our method\\
	\hline
    1 & $1024\times 2280$ & 1.324 sec & 6.203e-8 & 0.227 sec & 2.623e-9\\
	\hline
	2 & $4096\times 9120$ & 73.434 sec & 1.304e-6 & 3.233 sec & 1.577e-7\\
	\hline
	3 & $16384\times 36480$ & 4948.770 sec & 4.498e-5 & 76.473 sec & 2.361e-6\\
	\hline
\end{tabular}
\end{table}

\begin{figure}[h]
	\center{\includegraphics[scale=0.35]{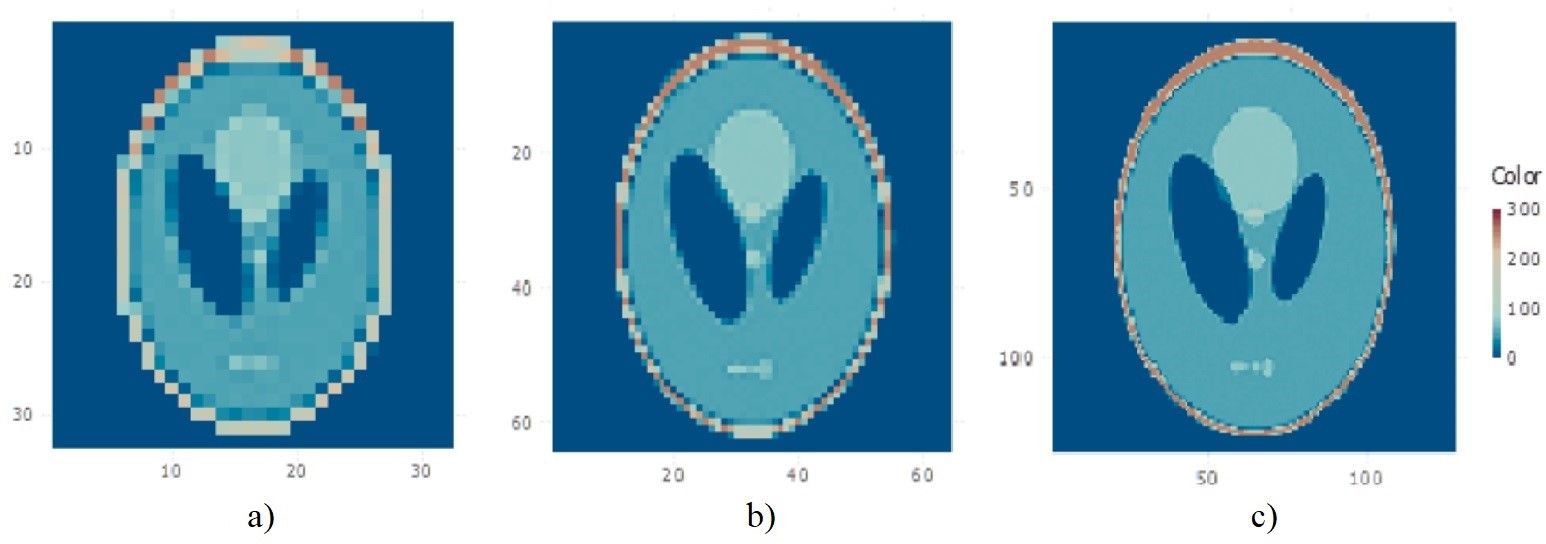}}
	\caption{The result of the reconstruction of the Shepp-Logan phantom. }
	\label{fig3}
\end{figure}

The proposed method is applicable for an arbitrary symmetric motion geometry of the source and detector, which ensures the symmetry of the reconstruction matrix. For example, in paper \cite{Xianguo} it was considered the problem of recovering images of nuclear waste in a cylindrical container using tomographic Gamma Scanning. Each reconstructed layer was divided into $72$ voxels using the polar coordinate system, and the voxel numbering was chosen in the order shown in Fig. 4a. If we redesign the voxel numbering by our proposed method, as shown in Fig. 4b, then the system of linear algebraic equations from \cite{Xianguo} can be divided into two independent systems with $36$ unknowns and the number of equations two times smaller than the original one.

\begin{figure}[h]
	\center{\includegraphics[scale=0.3]{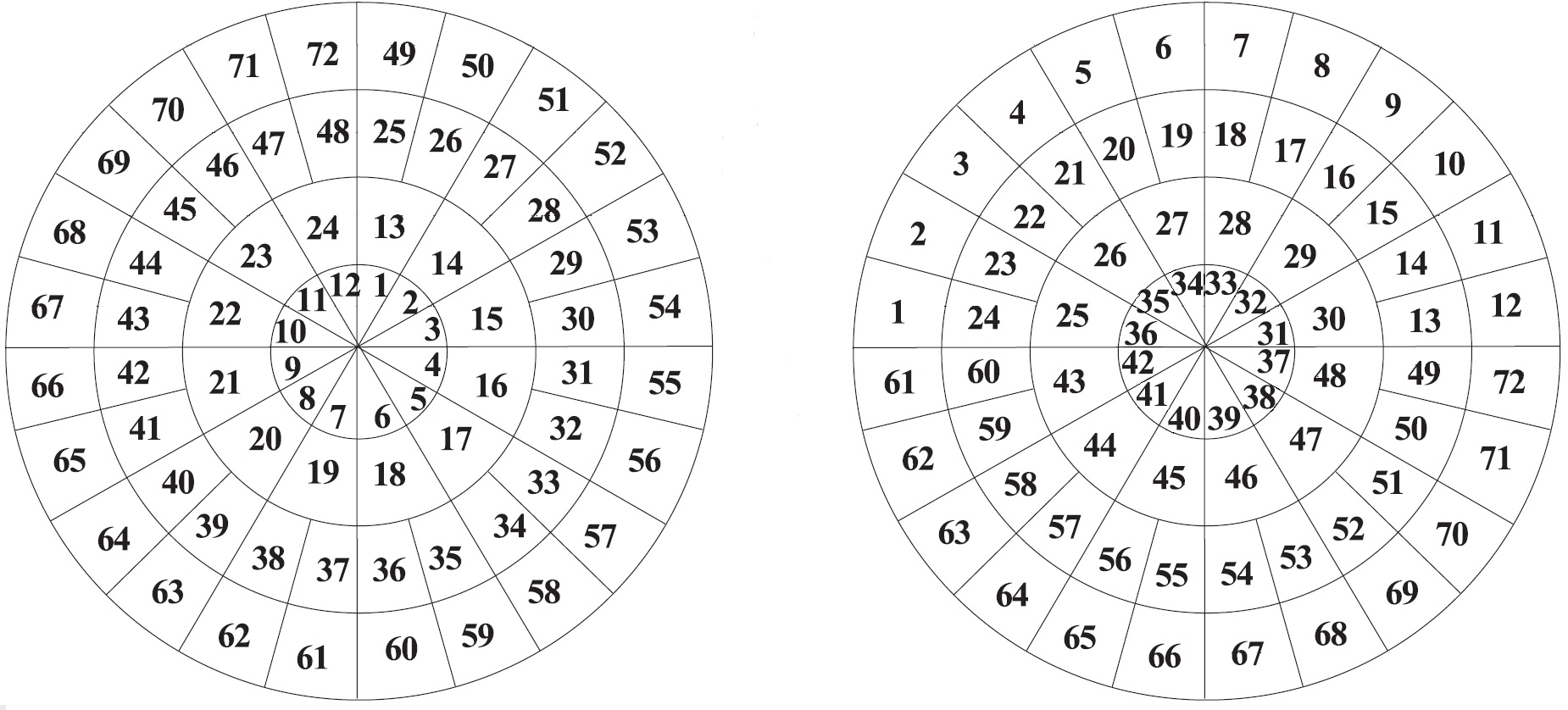}}
	\caption*{ a)   \hspace{75mm}     b)  }
	\caption{a) - voxel numbering proposed in \cite{Xianguo}, b) - voxel numbering according to our proposed method.}
	\label{fig4}
\end{figure}

\section{Conclusion}

A new method for solving systems of linear algebraic equations of a special type arising in solving problems of image reconstruction has been proposed. This method, due to a certain symmetry of the matrix and the choice of voxel numbering method for two-dimensional problems, allows us to divide the initial system of algebraic equations into two independent systems, which makes it possible to carry out the calculation in parallel. The size of the matrices of the resulting systems is 4 times less than the size of the original matrix, and these matrices are less dispersed.

The proposed method is applicable for an arbitrary symmetric motion geometry of the source and detector, which ensures the symmetry of the reconstruction matrix.

The combined use of this method and efficient algorithms, such as those proposed in \cite{006,ZhangSL,alikh, Xianguo}, for reconstructing a tomographic image will increase the speed of these algorithms at least 4 times, while maintaining the accuracy of the reconstruction.

This parallelization method can also be applied when solving multidimensional problems of mathematical physics considered in symmetric domains, the approximation of which leads to matrices symmetric in the sense of (\ref{ur2}).

\section{Data Availability}
The data used to support the findings of this study are available from the corresponding author upon request.

\section{Acknowledgements} 
The research results outlined in this paper were obtained with financial support from the Ministry of Education and Science of the Russian Federation, as part of the execution of the project entitled "Establishment and creation of hightech production for manufacturing of digital X-ray complex with the function of tomographic image synthesis", in accordance with the Government resolution of the Russian Federation \textnumero218 from 09.04.2010.


\begin{thebibliography}{5}

\bibitem{001} A.H. Andersen and A.C. Kak, Simultaneous algebraic reconstruction
technique (SART): A superior implementation of the art algorithm, Ultrason. Imag., vol. 6, no. 1, pp. 81--94, 1984, doi: 10.1016/0161-7346(84)90008-7

\bibitem{002} A. C. Kak, M. Slaney and G. Wang, Principles of computerized tomographic imaging, Med. Phys., vol. 29, no. 1, p. 107, 2002, doi: 10.1118/1.1455742

\bibitem{003} M. Jiang and G. Wang, Convergence of the simultaneous algebraic reconstruction technique (SART), IEEE Trans. Image Process. Publ. IEEE Signal Process. Soc., vol. 12, no. 8, pp. 957--961, Aug. 2003, doi: 10.1109/TIP.2003.815295

\bibitem{004} H. M. Hudson and R. S. Larkin, Accelerated image reconstruction using
ordered subsets of projection data, IEEE Trans. Med. Imag., vol. 13, no. 4, pp. 601--609, Dec. 1994, doi: 10.1109/42.363108

\bibitem{005} K. Mueller and R. Yagel, Rapid 3-D cone-beam reconstruction with
the simultaneous algebraic reconstruction technique (SART) using 2-D
texture mapping hardware, IEEE Trans. Med. Imag., vol. 19, no. 12,
pp. 1227--1237, Dec. 2000, doi: 10.1109/42.897815

\bibitem{006} B. Liu and L. Zeng, Parallel SART algorithm of linear scan cone-beam
CT for fixed pipeline, J. Xray Sci. Technol., vol. 17, no. 3, pp. 221--232,
2009, doi: 10.3233/XST-2009-0224

\bibitem{007} D. Xu, Y. Huang, and J. U. Kang, GPU-accelerated non-uniform fast
Fourier transform-based compressive sensing spectral domain optical
coherence tomography, Opt. Exp., vol. 22, no. 12, pp. 14871--14884,
2014, doi: 10.1364/OE.22.014871

\bibitem{008} X. Jia, B. Dong, Y. Lou, and S. B. Jiang, GPU-based iterative cone-beam
CT reconstruction using tight frame regularization, Phys. Med. Biol.,
vol. 56, no. 13, pp. 3787--3807, 2011, doi: 10.1088/0031-9155/56/13/004

\bibitem{009} P. B. Noel, A. M. Walczak, J. Xu, J. J. Corso, K. R. Hoffmann, and
S. Schafer, GPU-based cone beam computed tomography, Comput.
Methods Programs Biomed., vol. 98, no. 3, pp. 271--277, 2010, doi: 10.1016/j.cmpb.2009.08.006

\bibitem{010} W.-M. Pang, J. Qin, Y. Lu, Y. Xie, C.-K. Chui, and P.-A. Heng, Accelerating simultaneous algebraic reconstruction technique with motion compensation using CUDA-enabled GPU, Int. J. Comput. Assist. Radiol. Surg., vol. 6, no. 2, pp. 187--199, 2011, doi: 10.1007/s11548-010-0499-3

\bibitem{ZhangSL} SL. Zhang, GH. Geng, GH. Cao, YH. Zhang, BD. Liu, X. Dong, Fast Projection Algorithm for LIM-Based Simultaneous Algebraic Reconstruction Technique and Its Parallel Implementation on GPU, IEEE ACCESS, 2018, 6, 23007-23018, doi: 10.1109/ACCESS.2018.2829861

\bibitem{alikh} D.G. Kovtun, E.B. Gorbunova, A.A. Alikhanov, A.M. Apekov, A.O. Belyaev, A.S. Boldyreff, A.V. Ugolkov, Application of nonlinear voxel distribution grid for computational speed-up for linear tomosynthesis reconstruction, Proc. SPIE 10796, Electro-Optical Remote Sensing XII, 107960G (9 October 2018), doi: 10.1117/12.2513308


\bibitem{Xianguo} H. Aijing, T. Xianguo, S. Rui, Z. Honglong, An improved OSEM iterative reconstruction algorithm for transmission tomographic gamma scanning, Applied Radiation and Isotopes 142 (2018) 51-55.
doi: 10.1016/j.apradiso.2018.09.001

 

 



\end{thebibliography}
\end{document}